\documentclass[reqno]{amsart}
\usepackage{amsfonts,amscd, amssymb,  epsf, epsfig, graphicx}

\def\noi{\noindent}

\newcommand{\bea}{\begin{eqnarray*}}
\newcommand{\eea}{\end{eqnarray*}}

\newtheorem{lem}{LEMMA}[section]
\newtheorem{theo}[lem]{THEOREM}

\newtheorem{coro}[lem]{COROLLARY}
\newtheorem{prop}[lem]{PROPOSITION}

\newtheorem{remark}[lem]{Remark}
\newtheorem{main theorem}[lem]{MAIN THEOREM}

\newtheorem{main theo}[lem]{MAIN THEOREM}
\def\CQFD{\hfill \vrule width 7pt height 7pt depth 1pt}

\def\CC{{\rm\kern.24em\vrule width.02em height1.4ex depth-.05ex\kern-.26em C}}
\def\QQ{{\rm\kern.24em\vrule width.02em height1.4ex depth-.05ex\kern-.26em Q}}
\def\PP{{\rm\kern.24em\vrule width.02em height1.4ex depth-.05ex\kern-.26em P}}
\def\Rr{{\rm I\kern-.2em R}}
\def\ZZ{{\rm\kern.26em\vrule width.02em height0.5ex
depth0ex\kern.04em\vrule width.02em height1.47ex depth-1ex\kern-.34em Z}}
\def\BB{{\rm\kern.24em\vrule width.02em height1.4ex depth-.05ex\kern-.26em B}}
\def\RR{\hspace{.065in}\rm{\vrule width.02em height1.55ex
depth-.07ex\kern-.3165em R}}
\def\Ibb#1{{\rm I\kern-.23em#1}}

\def\noi{\noindent}

\def\noi{\noindent}

\newcommand {\be} {\begin{equation}}
\newcommand {\ee} {\end{equation}}

\newcommand {\beas} { \begin{eqnarray*}}
\newcommand {\eeas} {\end{eqnarray*}}

\begin{document}

\centerline{Density of Orbits in Complex Dynamics}

\bigskip

\centerline{J. E. Forn\ae ss, B. Stens\o nes}

\bigskip


\bigskip

\section{Introduction}

Let $F:\CC^m_0 \rightarrow \CC^m_0$ be a germ of a
holomorphic map with $F(0)=0, F'(0)={\mbox{Id}}.$
Dominique Cerveau ([C]) asked whether $F$ can have a dense
orbit. More precisely, does there exist a small ball
$\BB(0,\delta)=\BB^m_\delta=\{z\in \CC^m, \|z\|<\delta\}$ where $F$ is defined and a sequence
$S=\{p_n\}_{n=0}^\infty \subset \BB(0,\delta)$,
$F(p_n)=p_{n+1}\; \forall \; n$,  so that
$S \cap \BB(0,\epsilon)$ is dense in some smaller ball
$\BB(0,\epsilon), 0<\epsilon<\delta.$ 
Our main result is that this is not possible. We prove in fact a
more general result. Recall that a domain $\Omega$ 
in $\CC^{m}$is said
to have a Lipschitz boundary if we can locally write,
after a holomorphic rotation of coordinates,
$\Omega=\{{\mbox{Im}}z_{m}< r(z_1,\dots,z_{m-1},{\mbox{Re}}z_{m})\}$ where
$|r(x)-r(y)|\leq C \|x-y\|,$ the constant $C$ is called a Lipschitz constant
for $r.$ 

\begin{main theo}
Let $\Omega \subset \subset \CC^m$ be a domain with
a Lipschitz boundary. Suppose that $F:\Omega
\rightarrow \CC^m$ is a holomorphic map and let $V \subset \Omega$
be a nonempty open subset. Then $F$ does not have a dense orbit
in $V$, i.e. there does not exist a sequence 
$\{p_n\}_{n\geq 0}\subset \Omega$
for which $F(p_n)=p_{n+1}$ for all $n$ and
$\{p_n\}\cap V$ is dense in $V.$
\end{main theo}

The condition that $\Omega$ has Lipschitz boundary is not necessary,
see Remark 5.4. However, the following theorem shows that
the result is false if we have no boundary conditions.

\begin{theo}
There is a bounded domain $U$ in $\CC^m$ which is biholomorphic
to the unit ball and a holomorphic
map $F:U \rightarrow \CC^m$ with a dense orbit. So there is a dense sequence
$\{p_n\}_{n=0}^\infty \subset U$ with $F(p_n)=p_{n+1}$
for all $n.$
\end{theo}

The proofs depend on two results which are of independent
interest.

\begin{theo} Suppose that $\Omega \subset \subset \CC^m$
is a domain with Lipschitz boundary. Then if $p\in \partial
\Omega$ there are two not mutually exclusive possibilities:\\
\noi (i) The point $p$ is in the envelope of holomorphy of $\Omega.$\\
\noi (ii) The point $p$ is in a local peak set for plurisubharmonic
functions, i.e. there exists a neighborhood $U$ of $p$ and
a continuous function $\rho: U \cap \overline{\Omega}\rightarrow
\RR$ so that $\rho<0$ on $U \cap \Omega$ and $\rho$ is plurisubharmonic
on $U \cap \Omega$. Moreover, $\rho(p)=0.$
\end{theo}

\begin{theo}
Let $M^m$ be an $m-$dimensional
Hausdorff connected complex manifold with countable topology. Then
there exists a biholomorphic map $F: \BB^m(0,1) \rightarrow U$
from the unit ball onto an open subset $U \subset M$ so that
$M\setminus U$ has zero volume.
\end{theo}

We prove Theorem 1.4 in Section 2. The proof relies heavily on
previous joint work of the first author and E. L. Stout. They
proved, [FS1], that Theorem 1.4 is valid if one replaces the ball
with the polydisc of the same dimension. They also proved
covering theorems by balls, [FS2], in arbitrary complex manifolds. 
We would like to thank E. L. Stout for valuable
remarks during the work with Theorem 1.4.

\bigskip

Theorem 1.3 is proved in Section 3. The two main ingredients is
an estimate of the local envelope of holomorphy of a domain
with Lipshitz boundary together with the result by Demailly [D]
that Lipschitz pseudoconvex domains are hyperconvex, $i.e.$
that they allow a bounded plurisubharmonic exhaustion function.
The key result about envelopes of holomorphy is the following:

\begin{theo}
Let $r(z_1,\dots,z_{m-1},{\mbox{Re}}(z_{m}))$ be a real-valued
Lipschitz function on $\CC^{m-1}\times \RR$ with Lipschitz constant
$C$. Set
$$\Omega:=\{{\mbox{Im}}(z_{m})<r(z_1,\dots,z_{m-1},{\mbox{Re}}(z_{m}))\}.$$
Then the envelope of holomorphy $\tilde{\Omega}$ of $\Omega$
is schlicht and is of the form
$\tilde{\Omega}:=\{{\mbox{Im}}(z_{m})<\tilde{r}
(z_1,\dots,z_{m-1},{\mbox{Re}}(z_{m}))\},$
where $\tilde{r}\geq r$ is a Lipschitz function with the same
Lipshitz constant or $\tilde{r}\equiv \infty.$
\end{theo}

We prove Theorem 1.2 in Section 4. The proof goes by an explicit
construction of an expanding  map defined on a biholomorphic image of
the polydisc, and then using Theorem 1.4 to find a dense open subset
in this domain, biholomorphic to the ball. 

\bigskip

Finally we prove the Main Theorem in Section 5 by contradiction.
A key step is to use Theorem
1.3 to show that the iterates $F^n$ are well defined on $V$ with
image contained in $\overline{\Omega}.$ The key point is to show that
$F^n(V)$ is contained in $\tilde{\Omega}$ so that $F^{n+1}$ has a holomorphic 
extension.This allows
us to use normal families arguments.

\section{Proof of Theorem 1.4}

To prove Theorem 1.4 it suffices by Theorem 1 of [FS1] to assume that
$M=\Delta^m$. In fact we will prove the Theorem in the
slightly more general case, when $M=U$, a connected open set
in $\CC^m$. This is anyhow necessary as we change $U$ during the proof.

\begin{lem}
Let $U$ be a connected open set in $\CC^m$. Suppose
$F_j:\BB^m_{r_j} \rightarrow U, j=1,2$ are $1-1$ holomorphic maps
with $F_1(\BB^m_{r_1}) \cap F_2(\BB^m_{r_2}) =\emptyset.$
Let $0<s<1$. Then there exists a $1-1$ holomorphic map
$F_0:\BB^m_1\rightarrow U$ so that $F_0(\BB^m_1)
\supset F_1(\BB^m_{sr_1})\cup F_2(\BB^m_{\frac{sr_2}{32}})$
\end{lem}

\noindent{\bf Proof.} If $a$ is a real number, we shall understand by 
$\bf a$ the point
$(a,0,\ldots,0)$ in $\CC^m$. Also, for a positive number $t$, 
let $K_t$ be the union of the closed ball $\overline{\BB^m_1}$, 
the closed ball ${\mathbf 4}+\overline{\BB^m_t}$ and the straight 
line segment $[\mathbf 0,\mathbf 4]$. The set $K_t$ is thus a certain 
compact subset of $\CC^m$, which is polynomially convex when $t<3$. 
Fix a number $\sigma\in(s,1)$. Let $p\in F_1(b\BB^m_{\sigma r_1})$ 
be the point
$F_1({\mathbf{\sigma r_1}})$ 
and $q\in F_2(b\BB^m_{\sigma r_2})$ the point $F_2(-{\mathbf{\sigma r_2}})$. 
Let $\lambda$ be a smooth arc in $U$ that connects $p$ and $q$ and is 
otherwise disjoint from $F_1(\overline{\BB^m_{\sigma r_1}})
\cup F_2(\overline{\BB^m_{\sigma r_2}})$. Choose a smooth parametrization 
$f:[{\bf 1}, {\bf 2}]\rightarrow \lambda$ that takes $\bf 1$ to 
the point $p$ and $\bf 2$ to $q$.
The map $G:K_2\rightarrow \CC^m$ is defined by 
\begin{equation*}
G(z)=\begin{cases}
F_1(\sigma r_1 z)& \mbox{if $z\in\overline{\BB^m_1}$}\\ 
F_2(\frac{\sigma r_2 }{2}(z-{\mathbf{4}}))& 
\mbox{if $z\in{\mathbf{4}}+\overline{\BB^m_2}$}\\
f(x)& \mbox{if $z={\mathbf x}$ with $x\in[{\bf 1},{\bf 2}].$}\\
\end{cases}
\end{equation*}

With the proper choices of $\lambda$ and $f$, the function $G$ defined in 
this way is of class ${\mathcal {C}^\infty}$ on$[{\mathbf {0}}, 
{\mathbf{4}}]$. It is one-to-one on $K_2$.
The function $G$ can be approximated in the ${\mathcal{C}^1}$ sense on the 
set $K_2$ by polynomial maps $P:\CC^m\rightarrow\CC^m$ [FS1].
If $P$ is such a map that approximates $G$ well enough, then $P$ will be 
one-to-one on $K_2$ and can be chosen to be regular on a neighborhood of 
$K_2$. It will satisfy also the conditions

\noindent a) $P(K_2)\subset U$,

\noindent b) $P(K_2)\supset F_1(\BB^m_{sr_1})\cup F_2(\BB^m_{sr_2})$, and

\noindent c) $P\left({\mathbf{4}}+\BB^m_\frac{1}{16}\right)\supset 
F_2\left(\BB^m_{\frac{sr_2}{32}}\right)$.

Let $V\subset\CC^m$ be a neighborhood of $K_2$ on which $P$ is 1-1. Such 
neighborhoods exist, because $P$ is 1-1 on the compact set $K_2$ and is 
regular on a neighborhood of it. By Lemma 1 of [FS2], there is a domain 
$B\subset V$ that is a neighborhood of $K_\frac{1}{16}$ and that is 
biholomorphically equivalent to $\BB^m_1$.
 
The Lemma is proved.\\

\CQFD\\

\begin{coro}
Let $U$ be a connected open set in $\CC^m$. Suppose
$F_j:\BB^m_{r_j} \rightarrow U, j=1,2,\dots, N$ are finitely many
$1-1$ holomorphic maps
with $F_i(\BB^m_{r_i}) \cap F_j(\BB^m_{r_j}) =
\emptyset\; \forall\; i \neq j.$
Let $0<s<1$. Then there exists a $1-1$ holomorphic map
$F_0:\BB^m_1\rightarrow U$ so that $F_0(\BB^m_1)
\supset F_1(\BB^m_{sr_1})\cup_{j>1} F_j\left(\BB^m_{\frac{sr_{j}}{32}}\right)
$
\end{coro}

{\bf Proof of the Corollary:}
We prove the Corollary by induction. By the Lemma, the result
holds for $N=2.$ Assume the result holds for some $N>1.$
We will show that the result holds for $N+1$.
Let 
$F_j:\BB^m_{r_j} \rightarrow U, j=1,2,\dots, N+1$ be finitely many
$1-1$ holomorphic maps
with $F_i(\BB^m_{r_i}) \cap F_j(\BB^m_{r_j}) =
\emptyset\; \forall\; i \neq j.$
Let $0<s<1$. Pick a $\sigma, s<\sigma<1.$ Define $U_1:=
U \setminus \overline{F_{N+1}(\BB^m_{\sigma r_{N+1}})}.$ Then $U_1$
is a connected domain and 
$F_j:\BB^m_{r_j} \rightarrow U_1, j=1,2,\dots, N$  are finitely many
$1-1$ holomorphic maps
with $F_i(\BB^m_{r_i}) \cap F_j(\BB^m_{r_j}) =\emptyset\; \forall\; i \neq j
,i,j\leq N.$
Hence the inductive hypothesis applies to find a holomorphic map
$F_0:\BB^m_1 \rightarrow U_1$ so that
$F_0(\BB^m_1)
\supset F_1(\BB^m_{\sigma r_1})\cup_{1<j\leq N}
 F_j\left(\BB^m_{\frac{\sigma r_j}{32}}\right)$
It follows that the two maps
$F_0:\BB^m_1\rightarrow U$ and $F_{N+1}:\BB^m_{\sigma r_{N+1}}\rightarrow U$
have disjoint range.
Next, let $\lambda, s<\lambda<\sigma$ be arbitrary and pick 
$0<\tau<1$ large enough, then
$F_0(\BB^m_\tau)\supset 
 F_1(\BB^m_{\lambda r_1})\cup_{1<j\leq N}
 F_j\left(\BB^m_{\frac{\lambda r_{j}}{32}}\right).$

\medskip

Next we pick $\mu,0<\mu<1$ so large that 
$$F_0(\BB^m_{\mu\tau})\supset 
 F_1(\BB^m_{s r_1})\cup_{1<j\leq N} F_j\left(\BB^m_{\frac{ s r_{j}}{32}}
\right).$$
We apply the Lemma to the maps $F_0:\BB^m_\tau\rightarrow U$ and $F_{N+1}: 
\BB^m_{\sigma r_{N+1}}\rightarrow U.$ Pick $s', 0<s'<1$ so that $\mu<s'$
and $s<s' \sigma.$ Then there exists by Lemma 2.1 a $1-1$ holomorphic
map $F:\BB^m_1\rightarrow U$ so that

\bea
F(\BB^m_1) & \supset &  
F_0(\BB^m_{s'\tau})\cup F_{N+1}\left(\BB^m_{\frac{s'\sigma r_{N+1}}{32}}
\right)\\
& \supset & F_0(\BB^m_{\mu\tau})\cup 
F_{N+1}\left(\BB^m_{\frac{s r_{N+1}}{32}}\right)\\
& \supset &  
 F_1(\BB^m_{s r_1})\cup_{1<j\leq N} F_j\left(\BB^m_{\frac{ s r_{j}}{32}}
\right)
\cup 
F_{N+1}\left(\BB^m_{\frac{s r_{N+1}}{32}}\right)\\
\eea

\CQFD\\

Let $C$ denote a cube in $\CC^m$ and let $B$ denote a maximal
concentric ball in $C$. We let $\rho:=\rho_m:= \frac{|B|}{|C|}.$ Here
$|x|$ denotes the volume of $x$. Obviously $\rho_m$ is well-defined,
$i.e.$ independent of the size of $C$.

\begin{lem}
Let $V$ denote a bounded open set in $\CC^m.$ Then there exists a
finite number of disjoint balls $B_i \subset V$ so that
$|\cup B_i|> \frac{\rho|V|}{2}.$
\end{lem}

{\bf Proof:}
We can cover $\CC^m$ by a grid of cubes so that the volume
of the collection of cubes entirely contained in $V$ has
volume at least $\frac{|V|}{2}.$ Then pick the largest concentric ball
in each of these cubes.\\

\CQFD\\

\begin{lem}
Let $F:\BB^m_1\rightarrow U$ denote a $1-1$ holomorphic
map into a connected bounded open set $U\subset \CC^m$.
Let $0<s<1$ and $\epsilon>0.$ Then there exists a $1-1$
holomorphic map  
$F':\BB^m_1\rightarrow U$ so that 
$F'(\BB^m_1)\supset F(\BB^m_s)$ and $$|F'(\BB^m_1)|
\geq |F(B^m_1)|+\frac{\rho}{2*(32)^{2m}}\left(|U|-|F(\BB^m_1)|\right)
-\epsilon.$$
\end{lem}

{\bf Proof:}
If $r, 0<s<r<1$ is large enough,
then the volume of $F(\BB^m_1\setminus \BB^m_r)$ is
at most $\epsilon/2.$ 
The set $V:U \setminus F(\overline{\BB^m_r})$ is connected.
By Lemma 2.3 we can find finitely many disjoint
balls $B_j$ in $V$ with total volume
at least $\frac{\rho |V|}{2}.$ Applying Corollary 2.2
to these balls $B_j^m(p_j,r_j)$ and the ball $F(\BB^m_r)$ we find
a $1-1$ holomorphic image $F'(\BB^m_1)$ of the unit ball in $U$ which
contains $F(\BB^m_{\sigma r})$ and $B^m_j(p_j, \frac{\sigma r_j}{32})$
for any $\sigma, 0<\sigma<1, s<\sigma r$ we want. By choosing $\sigma$
large enough,
we can suppose that the volume of $F'(\BB^m_1)$ is at least
\bea
|F'(\BB^m_1)| & \geq & |F(\BB^m_r)|+ \frac{\rho |V|}{2*(32)^{2m}}-\epsilon/2\\
& \geq & 
|F(\BB^m_1)|+ \frac{\rho |V|}{2*(32)^{2m}}-\epsilon\\
& \geq & 
|F(\BB^m_1)|+ \frac{\rho (|U|-|F(\BB^m_1)|)}{2*(32)^{2m}})-\epsilon.\\
\eea

\noindent By the choice of $\sigma,$ we also have that
$F'(\BB^m_1) \supset F(\BB^m_{\sigma r}) \supset F(\BB^m_s).$\\

\CQFD\\

\begin{coro}
Let $U$ be a bounded connected open set in $\CC^m.$ Then there exists
a sequence of $1-1$ holomorphic maps $F_j:\BB^m_1\rightarrow \CC^m$
and numbers $0<s_j<1, j=1, \dots$ so that\\
\noindent (i) $F_{j+1}(\BB^m_{s_{j+1}}) \supset F_{j}(\BB^m_{s_{j}})\;
\forall\; j \geq 1$\\
\noindent (ii) $|F_{j+1}(\BB^m_{s_{j+1}})|
\geq |F_{j}(\BB^m_{s_{j}})|+\frac{\rho}{3*(32)^{2m}}\left(
|U|-|F_{j}(\BB^m_{s_{j}})|\right)\; \forall \; j \geq 1.$
\end{coro}

{\bf Proof:}
We prove the Corollary by induction. Let $F_1:\BB^m_1\rightarrow U$
be an arbitrary ball, set $s_1=1/2.$ Suppose we have found
$F_1,\dots,F_j$ and $s_1,\dots,s_j.$ We apply Lemma 2.4 to
the map $F=F_j$ and $s=\frac{s_j+1}{2}>s_j$. Let
$\epsilon= \frac{\rho}{12*(32)^{2m}}\left(|U|-|F_j(\BB^m
_{s_j})|\right).$
Then pick $F'$ as in Lemma 2.4 and set $F_{j+1}=F'$ Then
$$F_{j+1}(\BB^m_1)\supset F_j(\BB^m_s) \supset \supset F_j(\BB^m_{s_j})\;
{\mbox{and}}$$

\bea
|F_{j+1}(\BB^m_1)| & \geq & |F_j(\BB^m_1)|+
\frac{\rho}{2*(32)^{2m}}\left(|U|-|F_j(\BB^m_1)|\right)-\epsilon\\
& = &  |F_j(\BB^m_1)|(1-
\frac{\rho}{2*(32)^{2m}})+\frac{\rho}{2*(32)^{2m}}|U|-\epsilon\\
& \geq & 
|F_j(\BB^m_{s_j})|(1-
\frac{\rho}{2*(32)^{2m}})+\frac{\rho}{2*(32)^{2m}}|U|-\epsilon\\
& = & 
|F_j(\BB^m_{s_j})|+
\frac{\rho}{2*(32)^{2m}}\left(|U|-|F_j(\BB^m_{s_j})|\right)-\epsilon\\
& = & 
|F_j(\BB^m_{s_j})|+
\frac{5\rho}{(12)*(32)^{2m}}\left(|U|-|F_j(\BB^m_{s_j})|\right)\\
& > & |F_j(\BB^m_{s_j})|+
\frac{\rho}{3*(32)^{2m}}\left(|U|-|F_j(\BB^m_{s_j})|\right)\\
\eea

Finally we choose $0<s_{j+1}<1$ large enough that

$$F_{j+1}(\BB^m_{s_{j+1}}) \supset F_j(\BB^m_{s_j})\;
{\mbox{and}}$$

$$|F_{j+1}(\BB^m_{s_{j+1}})|\geq |F_j(\BB^m_{s_j})|+
\frac{\rho}{3*(32)^{2m}}\left(|U|-|F_j(\BB^m_{s_j})|\right)$$

\CQFD\\

\begin{prop}
Let $U$ be a bounded connected open set in $\CC^m.$ Then there exists
a sequence of $1-1$ holomorphic maps $F_j:\BB^m_1\rightarrow \CC^m$
and numbers $0<s_j<1, j=1, \dots$ so that\\
\noindent (i) $F_{j+1}(\BB^m_{s_{j+1}}) \supset F_{j}(\BB^m_{s_{j}})\;
\forall\; j \geq 1$\\
\noindent (ii) $|F_{j}(\BB^m_{s_{j}})| \rightarrow |U|$
\end{prop}

{\bf Proof:}
Let $F_j$ and $s_j$ be as in the Corollary 2.5.
Set $A:=\lim_j |F_{j}(\BB^m_{s_{j}})|.$ From (ii) it follows
that 

$$
|U| \geq A \geq A+\frac{\rho}{3*(32)^{2m}}(|U|-A).
$$

This is impossible if $A<|U|.$\\

\CQFD\\

{\bf Completion of the Proof of Theorem 1.4:}
The union $\cup_{j} F_j(B^m_{s_j})$ is biholomorphic to the unit ball.
This follows from [FS1] since the union is contained in a
Kobayashi hyperbolic set.\\

\CQFD\\

\section{Proof of Theorem 1.3}

We prove first Theorem 1.5.

{\bf Proof:}
Let $c$ be a Lipschitz constant of $\phi.$ Then

$$
\frac{|\phi(v,z)-\phi(v_0,z_0)|}{\|(v,z)-(v_0,z_0)\|}
\leq c.$$

 \noi This is equivalent to the following:
Let $\Gamma=\{u=\phi(v,z)\}$ and choose a point $(u_0+iv_0,z_0)\in \Gamma,$
then $\Gamma$ does not intersect the cone 

$$
{\mathcal C}_{(z_0,u_0,v_0)}:=\{(z,u+iv);|u-u_0|> c \|(z,v)-(z_0,v_0)\|.
$$

\noi in the point $(u_0+iv_0,z_0).$ Moreover let $(\alpha,\beta)
\in \CC^{m-1}\times \RR$ be such that
$\|(\alpha,\beta)\|<\frac{1}{c}$, then all points
$(z_0,u_0+iv_0)-t(\alpha,1+i\beta)$ where $0<t$ are contained
in $\Omega.$ Let

$$
 H_{t,\alpha,\beta}(z,u+iv)=(z,u+iv)-t(\alpha,1+i\beta).
$$

\noi For fixed $t,\alpha$ and $\beta$ this is a holomorphic map and

$$
(*)\;\;\;\;\;\;\;\;\;\;\;\;\; H_{t,\alpha,\beta}(\Omega) \subset \Omega.
$$

\bigskip

We want to show that $\tilde{\Omega}$ also has this property $(*)$.
Since $\Omega$ is bounded by a graph it follows that $\tilde{\Omega}$
is schlicht and also is bounded by a graph (see [S]). In fact,
$\tilde{\Omega}=\cup_{j=0}^\infty \Omega_j, \Omega_0=\Omega,$ where $\Omega_j$
consists of all the points that can be reached from $\Omega_{j-1}$
by pushing discs with boundaries in $\Omega_{j-1}$. We make this more
precise: Let $H$ denote the Hartogs figure

$$
H:=\{p\in 
\overline{\Delta}_0 \cup \cup_{0\leq s\leq 1}\partial \Delta_s\},\;
\Delta_s=\{(z,t)\in \CC \times \RR; t=s, |z|\leq 1\}\subset \CC^2.
$$

Moreover, $\hat{H}$ denotes the convex hull of $H.$ Let $\Phi$ denote
any rank $2$, one-to-one holomorphic map defined on some (not fixed)
neighborhood of $\hat{H}$ in $\CC^2$ with values in $\CC^m.$ Then

$$
\Omega_j:= \cup_{\{\Phi; \Phi(H)\subset \Omega_{j-1}\}}
\Phi(\hat{H}).
$$

\medskip

Now we only need to show that if $\Omega_{j-1}$ has property $(*)$
so does $\Omega_{j}.$ 

\medskip

Observe that if $\Phi(H)\subset \Omega_{j-1}$, then
$H_{t,\alpha,\beta}(\Phi(H))= (H_{t,\alpha,\beta}\circ \Phi)(H))
\subset \Omega_{j-1}.$ Hence if $p\in \Omega_j,$
then $H_{t,\alpha,\beta}(p)\in \Omega_j.$

\bigskip

Now $\partial \tilde\Omega =\tilde{\Gamma}$ where

$$
\tilde{\Gamma}=\{(z,u+iv); u=\tilde{\phi}(z,v)\}.
$$

Obviously, by $(*)$, if $\tilde{\phi}(z,v)$ is infinite for
some $(z,v),$ then $\tilde{\phi} \equiv \infty$ and we are done.
We assume hence that $\tilde{\phi}$ has values in $\RR$ only.

We need to show that if $(z_0,u_0+iv_0)\in \tilde{\Gamma},$ then
$\tilde{\Gamma}$ does not intersect the cone ${\mathcal C}_{(z_0,u_0,v_0)}.$

We know that the cone does not meet $\tilde{\Gamma}$
in points $(z,u+iv)$ where $u<u_0$ since these are 
interior points of $\tilde{\Omega}.$

\medskip

Assume that $(z_1,u_1+iv_1)\in \tilde{\Gamma}\cap {\mathcal C}_{(z_0,u_0,v_0)}$
and $u_1>u_0.$ We want to show that this implies that
$(z_0,u_0+iv_0)\in \tilde{\Omega}$ which contradicts the fact that this is a
boundary point. But the fact that $(z_1,u_1+iv_1)\in \partial \tilde{\Omega}$
implies that the set

$$
\{(z,u+iv)=(z_1,u_1+iv_1)-t(1,\alpha,\beta); \|(\alpha,\beta)\|
<\frac{1}{c},0<t<1\}
$$

\noi is contained in $\tilde{\Omega}$ and a simple computation
shows that $(z_0,u_0+iv_0)$ is an interior point of this set.\\

\CQFD\\

{\bf Proof of Theorem 1.3:}
Pick a point $p\in \partial \Omega$. If $p$ is an interior point
of the envelope of holomorphy then we are done. So suppose not.
Since $ \Omega$ has Lipschitz boundary, we can rotate and translate coordinates
so that $p=0$ and $\Omega$ is locally given by
$\Omega=\{(z,u+iv);u<\phi(z,v),\|(z,v)\|<\rho\}$ for some Lipshitz function $\phi.$
let $A>\phi(0)$ and pick $L>0$ large enough that
$A-L\|(z,v)\|<\phi(z,v)$ if
$\rho/2<\|(z,v)\|<\rho.$ We globalize $\Omega$ by setting

\bea 
\Omega'  & = & \{(z,u+iv)\in \CC^m; u-L\|(z,v)\|\; {\mbox{if}}\;
\|(z,v)\|>\rho/2,\\
& & 
u<\min\{A-L\|(z,v)\|,\phi(z,v)\}\;{\mbox{if}}\;
\|(z,v)\|\leq \rho/2\}.
\eea

Letting $\tilde{\Omega}'$ be the envelope of holomorphy
of $\Omega'$, $\tilde{\Omega}'=\{u<\tilde{\phi}(z,v)\},$
we get $\tilde{\phi}(0)=\phi(0)=0,
\tilde{\phi}\geq \phi$ near $0.$
By a
theorem by Demailly ([De]) it follows that there is a
locally defined continuous function $\rho$ on the set
$\{u\leq \tilde{\phi}(z,v)\}$ which is plurisubharmonic on
 $\{u< \phi(z,v)\} \subset\{u< \tilde{\phi}(z,v)\}$
and with $\rho(0)=0.$\\

\CQFD\\

\section{Proof of Theorem 1.2}  

We prove the theorem in steps.

\bigskip

Step 1:\\

Let $R'$ be the closed rectangle in $\CC$ with corners
$$(0,0), (2\sqrt{2}\pi,0), (0,4\pi), (2\sqrt{2}\pi,4\pi).$$ This is covered
by two closed rectangles $R_1',R_2'$ where $R_1'$ has corners
$$(0,0), (2\sqrt{2}\pi,0), (0,2\pi), (2\sqrt{2}\pi,2\pi)$$
\noindent and $R'_2$
has corners  $(0,2\pi), (2\sqrt{2}\pi,2\pi), (0,4\pi), (2\sqrt{2}\pi,4\pi).$

Let $F_1:R_1'\rightarrow R', F_1(z)=\sqrt{2}i z+2\sqrt{2}\pi$
and $F_2:R_2'\rightarrow R', F_2(z)=\sqrt{2}iz+4\sqrt{2}\pi.$ Then
the $F_j$ are homeomorphisms. 
Next let $R$, $R_j$ denote the interiors of the same rectangles.
Let $\Omega:=R_1\cup R_2$ and set $F:\Omega \rightarrow R, F(z)=
F_1(z), z\in R_1, F(z)=F_2(z), z\in R_2.$ 
Let $p$ be any point in $\Omega$ and let $U\subset \Omega$ be some
neighborhood of $p.$ Then there exists a connected
open subset $V\subset U$ and an
integer $N\geq 1$ so that $F^j(V)\subset \Omega$ for each $1\leq j
\leq N,F^N(V)=R_1$ or $R_2$ and hence $\Omega \subset F^{N+1}(V).$
From this observation it is easy to inductively find a dense
orbit $\{p_n\}$ in $\Omega.$ Moreover, we can choose $\{p_n\}$
to be contained in any open subset of full measure or just dense.\\

Step 2:\\

Let $D$ denote a disc,

$$
D:=\{w\in \CC; 1<|w|<e^{2\sqrt{2}\pi},w\notin \RR^+\setminus {0}\}.
$$

Then the maps $\tau_j:R_j \rightarrow D,
\tau_j(z)=e^z$ are biholomorphisms, $j=1,2.$
Suppose that $z_j\in R_j, z_2=z_1+2\pi i.$ Then $\tau_1(z_1)=
\tau_2(z_2).$
Notice also that
\bea
F_2(z_2) & = & \sqrt{2}iz_2+4\sqrt{2}\pi\\
& = & \sqrt{2}iz_1-2\sqrt{2}\pi+4\sqrt{2}\pi\\
& = & F_1(z_1)\\
\eea

Hence, we can define $\Phi:D \rightarrow \CC$ by $\Phi(w)=e^{F_j(z_j)},
j=1\; \mbox{or}\; 2,w=e^{z_j},z_j\in R_j.$
Then $\Phi$ is a well defined holomorphic map. Let $\{p_n\}_{n=0}^\infty$
be a dense orbit of $F$ in $\Omega.$ Define $q_n:=e^{p_n}\in D.$
Then $\Phi(q_n)=e^{F_j(p_n)}$ if $p_n\in R_j.$
Hence $\Phi(q_n)=e^{F(p_n)}=e^{p_{n+1}}=q_{n+1}.$ Hence
$\{q_n\}\subset D$ and is a dense orbit for $\Phi.$ \\

Step 3: \\

Let $D^m:= D \times \cdots \times D$, $m$ times and define
$\Psi:D^m \rightarrow \CC^m$ by
$\Psi(z_1,\dots,z_m)= (\Phi(z_1),\dots, \Phi(z_m)).$
Then we can find a dense sequence $\{p_n\}\subset D^m$ contained
in any given dense open subset.\\

Step 4:\\

To finish the proof of Theorem 1.2 we use Theorem 1.4 to find
an open subset of $D^m$ of full measure and biholomorphic
to the unit ball.\\

\CQFD\\

\section{Proof of the Main Theorem}

To prove the Theorem, we assume to the contrary that there are such $V$ and
$\{p_n\}.$ We can assume $V$ is connected and $(0)\in V.$
We let $\tilde{\Omega}$ denote the (possibly multisheeted)
envelope of holomorphy of $\Omega$ and $\Omega^*\subset
\tilde{\Omega}$ consists of $\Omega \cup \{q\in \partial \Omega;
q\in \tilde{\Omega}\}.$ The map $F$ extends holomorphically
to a map $\tilde{F}:\tilde{\Omega} \rightarrow \CC^m.$ We denote by
$F^*$ the restriction of $\tilde{F}$ to $\Omega^*.$ So $F^*$ is a continuous
extension of $F$ to $\Omega^*$ and in turn $F^*$ extends holomorphically.
We define $\Omega_k:=\{z=z_0\in \Omega; z_1:=F^*(z)\in \Omega^*,
z_2:=F^*(z_1)\in \Omega^*,\dots, z_k:=F^*(z_{k-1})\in \Omega^*\}.$ 
Note that in general the $\Omega_k$ are not open sets but
the maps $(\tilde{F})^{k+1}$ are well defined holomorphic maps in small
neighborhoods of the $\Omega_k.$

\begin{lem}
We have that $V \subset \Omega_k$ for all $k.$
\end{lem}

{\bf Proof of the Lemma:}
We show first that $F(V) \subset \Omega_1.$ 
If $p_k\in V$, then $F(p_k)=p_{k+1}\in \Omega.$ Since
$\{p_n\}\cap V$ is dense in $V$, it follows by continuity
that $F(V) \subset \overline{\Omega}.$ Suppose there is a point
$q=F(p)\in \partial \Omega,p\in V.$ We will show that $q\in \Omega^*.$
If $q \notin\Omega^*$, we can find, by Theorem 1.3, an open
neighborhood $U(q)$ and a continuous function $\rho :U\cap \overline{\Omega}
\rightarrow (-\infty,0]$ which is plurisubharmonic on $U \cap \Omega$,
$\rho<0$ on $U \cap \Omega$ and $\rho(q)=0.$ Since $\Omega$
has Lipschitz boundary, we can find an outward vector $\bf{n}$ in a possibly
smaller neighborhood of $q.$ 
Let $\Phi_\epsilon, \epsilon>0$ denote the inward translation
by $\epsilon$ in the direction of $\bf n.$ Then there exists
a fixed neighborhood $q\in W\subset V$ so that
$\rho\circ \Phi_\epsilon\circ F$ is plurisubharmonic on
$W$ for all small enough $\epsilon$ and moreover
$\rho \circ \Phi_\epsilon \circ F \rightarrow \rho\circ F$
uniformly as $\epsilon \searrow 0.$ It follows that
$\rho\circ F$ is plurisubharmonic on $W$ and since this
function takes its maximum at $q$ it must be constant on $W.$
This is a contradiction since $\rho\circ F<0$ at each point
of $\{p_n\}\cap U.$ This shows that $V \subset \Omega_1.$ \\

Suppose we have
shown that $V \subset \Omega_k.$ Then $F^{k+1}$ extends holomorphically
to $V$. The same argument as above shows that the image
of the extension is contained in $\overline{\Omega}.$ Furthermore,
applying theorem 1.3 again in the same way as above, we get that
the image is actually contained in $\Omega^*.$ This shows that
$V \subset \Omega_{k+1}$, completing the induction.\\

\CQFD\\

{\bf Proof of the Theorem:}\\

Since $(0)\in V,$ we can
choose sequences $m(k)<\ell(k)$ so that $p_{m(k)}, p_{\ell(k)}
\rightarrow 0.$ We set $n(k)=\ell(k)-m(k)$. We can assume that
$n(k) \nearrow \infty$. By the above Lemma, the map $F^{n(k)}$ has a
holomorphic extension to $V$ and $F^{n(k)}(p_{m(k)})=p_{\ell(k)}.$ 
Restricting to a thinner subsequence if necessary,
we can assume that $F^{n(k)} \rightarrow G: V \rightarrow \CC^m,$
$G(0)=0.$ 

\medskip

For any positive integer $\tau$ the sequence
$F^{\tau n(k)}=({F^{n(k)}})^{\tau}$ converges in a small
enough neighborhood to $G^\tau$ and by a normal families
argument, this convergence extends to all of $V$ and a suitable
extension of $G^\tau$ to $V$. The image of $V$ is still
contained in $\overline{\Omega}.$

\bigskip

We get:

\begin{lem}
Suppose that $F^{n_k}:V \rightarrow \overline{\Omega}$ is any sequence of
iterates converging to some map $G:V \rightarrow \Omega, G(0)=0,$
then all iterates $G^\tau$ can also be analytically continued
from $0$ to a map from $V$ to $\overline{\Omega}$. Moreover, if
$G^{\tau_s}$ is any convergent subsequence to a map
$H: V\rightarrow \CC^m$, then there is a subsequence
$F^{m(s)}$ converging to $H$ on $V$, moreover,
$H(0)=0, H(V) \subset \overline{\Omega}.$
\end{lem}

Using this Lemma we can refine the properties of
a limit $G$ of some subsequence $F^{n_k}$.\\

Let $\{\lambda_j\}_{j=1}^m$ denote the eigenvalues
of $G'(0)$. Since $G^\tau(V) \subset \overline{\Omega}$ for all
$\tau$, it follows that all eigenvalues have modulus $\leq 1.$
We can assume that
$1\geq |\lambda_1| \geq \cdots \geq |\lambda_m|.$
Using the Lemma again, we can assume that $\lambda_i=1,
i\leq j, \lambda_i=0, i>j$ for some $0 \leq j\leq m.$ 
If we put the matrix $G'(0)$ in Jordan normal form, we see that
for the blocks with eigenvalue $1$ there can be no off diagonal term
because these would grow under iteration violating Schwarz's Lemma.
Also if we look at the off diagonal terms for blocks with eigenvalue
zero, these vanish after iteration, so hence we can assume that there are
no off-diagonal terms.

\bigskip

Case (i): $\lambda_1=\cdots =\lambda_m=0.$\\

Then there exists $0<s<r$ and $k$ sufficiently large that
$F^{n_k}(\BB(0,r))\subset \BB(0,s)\subset \BB(0,r)\subset V$ and
$\|(F^{n_k})'(a)(v)\|\leq c\|v\|$ for some $c<1$ and
for all $a\in \BB(0,r).$
Hence $F^{n_k}$ is a contraction and hence has an attracting fixed
point. Therefore $F$ has an attracting periodic orbit.
This is impossible since there is a dense orbit.\\

Case (ii): $\lambda_1= \cdots =\lambda_m=1.$\\

In this case $G'(0)=\mbox{Id}$. If $G$ has nonzero higher order
terms in the Taylor expansion around the origin, the iterates of $G$
cannot be a normal family on $V$, so $G(z)\equiv z$ on $V$.

\medskip

Let $\hat{V}:=V \cup \cup_{n\geq 1}{\tilde{F}}^n(V).$ Then $\hat{V}$
is open and $\tilde{F}(\hat{V})\subset \hat{V}.$

\begin{lem}
The set $\hat{V}$ has finitely many connected components
$U_1,\dots,U_r=U_1$ with $\tilde{F}(U_i)\subset U_{i+1}.$ In fact
the map $\tilde{F}:U_i\rightarrow U_{i+1}$ is a biholomorphism and
${\tilde{F}}^{n_k} \rightarrow {\mbox{Id}}$ on each $U_i.$
\end{lem}

{\bf Proof of the Lemma:}
Clearly $\tilde{F}^m$ must be $1-1$ on $V$ for each iterate $\tilde{F}^m$.
Hence the Jacobians of each iterate is everywhere nonzero and
each iterate is an open map. Hence $\hat{V}$ is an open set.
Moreover, $\tilde{F}(\hat{V})\subset \hat{V}$. In particular,
each connected component is mapped into another connected component
by $\tilde{F}$. Let $U_1$ be the connected component containing $V$ and
define inductively connected components $U_i$ by $\tilde{F}
(U_i)\subset U_{i+1}$.
Since $\tilde{F}^{n_k}\rightarrow {\mbox{Id}}$ on $V$ it follows that for a minimal
$r, U_r=U_1.$ Hence there are only finitely many connected components.
The condition that $\tilde{F}^{n_k}\rightarrow {\mbox{Id}}$ on $V$ implies
that $F:\hat{V} \rightarrow \hat{V}$ is onto.
Hence each map $\tilde{F}:U_i\rightarrow U_{i+1}$ is onto. 
By the identity theorem
it follows that $\tilde{F}^{n_k}\rightarrow {\mbox{Id}} $ on each $U_i$. Hence
 $\tilde{F}:U_i\rightarrow U_{i+1}$ is $1-1$ and hence a biholomorphism.\\

\CQFD\\

This implies that ${\tilde{F}}^r:U_1 \rightarrow U_1$ is an automorphism
with a dense orbit of a bounded domain. The automorphism group is a
Lie group and the closure of the set of iterates of ${\tilde{F}}^r$ 
is a commutative
closed subgroup, isomorphic to a product of a torus and Euclidean
space and a finite commutative group. 
Since there is a sequence of iterates of ${\tilde{F}}^r$
converging to the identity, the Euclidean component must vanish, 
but then the orbit
of any point is nowhere dense, a contradiction.\\

Case (iii): $\lambda_1= \cdots =\lambda_j=1, \lambda_{j+1}=\cdots =
\lambda_m=0$
for some $1 \leq j <m.$
Then $G$ is of the form
$$
(**)\;\;\;\;\;\;\;
G(z_1,\dots,z_m)=(z_1+f_1,\dots,z_j+f_j,f_{j+1},\dots,f_m), f_i=
{\mathcal O}(\|z\|^2)$$

Fix a small neighborhood $W, 0\in W \subset \subset U_1$.
Let $U \subset \subset W$ and suppose that ${\tilde{F}}^r(U) \subset W.$ 
Then we have the 
functional equation $${\tilde{F}}^r\circ G(z)=G \circ {\tilde{F}}^r(z), 
z\in W.$$

Let $z_\delta:=(0,\c\dots,0,\delta)$ for small enough $\delta \neq 0.$
Using the density of orbits and equicontinuity there exists for any 
$\epsilon>0$ a neighborhood $U(0)\subset W$ and an integer $r_\epsilon>1$
so that $\|{\tilde{F}}^{r_\epsilon}(0)-z_\delta\|\leq \epsilon$ and 
${\tilde{F}}^r(U)\subset W.$
Since $G(0)=0$, ${\tilde{F}}^{r_\epsilon}(0)=G({\tilde{F}}^{r_\epsilon}(0))$. 
This is 
clearly impossible by (**). This completes the proof.\\

\CQFD\\

\begin{remark}
We can relax the Lipschitz condition of The Main Theorem if we add more
hypotheses on $F:$ Let $\Omega$ be any domain
with $\Omega={\mbox{int}}(\overline{\Omega})$
and assume in addition that the map $F$ has discrete fibers. 
Then $F$ has no dense orbit.
\end{remark}

\bigskip

\noi John Erik Forn\ae ss, Berit Stens\o nes\\
Mathematics Department\\
East Hall\\
University of Michigan\\
Ann Arbor, Michigan 48109\\
USA\\
fornaess@umich.edu, berit@umich.edu\\

\end{document}